\documentclass{amsart}

\usepackage{latexcad}

\newtheorem{theorem}{Theorem}[section]
\newtheorem{lemma}[theorem]{Lemma}
\newtheorem{proposition}[theorem]{Proposition}
\newtheorem{corollary}[theorem]{Corollary}
\theoremstyle{definition}
\newtheorem{definition}[theorem]{Definition}

\begin{document}

\title[Distances of Groups of Prime Order]{Distances of Groups of Prime Order}

\author{Petr Vojt\v echovsk\'y}

\address{Department of Algebra, Faculty of Mathematics and Physics,
Charles University, Sokolovsk\'a 83,
Prague, Czech Republic}

\curraddr{Department of Mathematics, Iowa State University, Ames, IA, U.S.A.}
\email{petr@iastate.edu}

\thanks{
 While working on this paper the author has been partially
 supported by the University \linebreak Development Fund of Czech Republic,
 grant number 1379/1998
}

\def\oa{\circ}
\def\ob{\ast}
\def\iso{\simeq}
\def\to{\longrightarrow}
\def\defining#1{{\it #1}}

\maketitle

\section{Introduction}

Let $G$ be a~finite set with $n$ elements, and $G(\oa)$, $G(\ob)$ two groups
defined on $G$. Their (\defining{Hamming}) \defining{distance} is the number of
pairs $(a,b)\in G\times G$ for which $a\oa b\ne a\ob b$. Let us denote this value
by $dist(G(\oa),G(\ob))$.

It is not difficult to show that
$dist(\underbar{\phantom{g}},\underbar{\phantom{g}})$ is a~metric on the set of
all groups defined on $G$. In fact, when $G_n$, $G_m$ are two groups of
different orders $n$ and $m$, respectively, and $dist(G_n,G_m)$ is defined
simply by $max\{n^2,m^2\}$, then
$dist(\underbar{\phantom{g}},\underbar{\phantom{g}})$ is a~metric on all finite
groups (defined on some fixed sets).

Similar ideas were first introduced by L.~Fuchs in \cite{Fu}. He asked about
the maximal number of elements, which can be deleted at random from a~group
multiplication table $M$, so that the rest of $M$ determines $M$ up to
isomorphism, or even allows a complete reconstruction of $M$. These two numbers
have been denoted by $k_1(M)$ and $k_2(M)$.

J.~D\'enes shows in \cite{Den1} that $k_2(M)=2n-1$, not including abelian
groups of order $4$ and $6$. His proof (published also in \cite{Den2}) was
fixed by S.~Frische in \cite{Fri}. She also found correct values of $k_2(M)$
for abelian groups of order $4$ and $6$ --- these are equal to $3$ and $7$.
Surprisingly, the value of $k_2(M)$ does not depend on structure of $M$ at all.

\begin{definition}
Let $G(\oa)$ be a~group. Then
\begin{displaymath}
    \delta(G(\oa))=min\{dist(G(\oa),G(\ob)); G(\ob)\ne G(\oa)\}
\end{displaymath}
is called \defining{Cayley stability of $G(\oa)$}. In similar manner, put
\begin{align*}
    \mu(G(\oa))&=min\{dist(G(\oa),G(\ob)); G(\ob)\iso G(\oa)\ne G(\ob)\},\\
    \nu(G(\oa))&=min\{dist(G(\oa),G(\ob)); G(\ob)\not\iso G(\oa)\},
\end{align*}
and call these numbers \defining{Cayley stability of $G(\oa)$ among isomorphic
groups}, \defining{Cayley stability of $G(\oa)$ among non-isomorphic groups},
respectively. Note that $\nu(G(\oa))$ is defined only when $n$ is not a~prime.
\end{definition}

\begin{definition}
Let $f:H\to K$ be a~mapping between two groups $H$, $K$.
\defining{Distance of $f$ from a~homomorphism} is the number $m_f$ of pairs
$(a,b)\in H\times H$ at which $f$ does not behave as a~homomorphism, i.e.
$f(ab)\ne f(a)f(b)$.
\end{definition}

When both operations $\oa$ and $\ob$ are fixed, and $g$ is an element of $G$, we
shall use $d(g)$ to denote the cardinality of $\{h\in G; g\oa h\ne g\ob h\}$.

\section{Some known facts}

Relatively few facts are known about $\nu(G(\oa))$. One can prove that
$v(E_{2^n})=2^{2n-2}$, where $E_{2^n}$ is the elementary abelian $2$-group of order
$2^n$ (see \cite{DrP}).
More generally, when $G(\oa)$, $G(\ob)$ are two groups of order $n$ with
$d(G(\oa),G(\ob))<n^2/4$, then their Sylow $2$-subgroups must be isomorphic
(see \cite{Dr2}).

The Cayley stability is known for any group $G(\oa)$ of order $n\ge 51$
(main result of \cite{DrM}), and is
equal to $\delta_0(G(\oa))$, where, using words of \cite{3W},
\begin{displaymath}
    \delta_0(G(\oa))=\begin{cases}
        6n-18& \text{if $n$ is odd},\\
        6n-20& \text{if $G(\oa)$ is dihedral of twice odd order},\\
        6n-24& \text{otherwise}.
    \end{cases}
\end{displaymath}
Cayley stability of $G(\oa)$ is less than or equal to $\delta_0(G(\oa))$
whenever $n\ge 5$ (for more details see \ref{drapal2}). Moreover, the nearest
group $G(\ob)$ must be isomorphic to $G(\oa)$. As \ref{drapal2} says, when
$f:G(\oa)\to G(\ob)$ is an isomorphism, then $f$ is a~transposition. This means
that $\mu(G(\oa))<\nu(G(\oa))$ holds for all groups of order at least $51$.
However, $\mu(G(\oa))<\nu(G(\oa))$ is not true in general; the exceptions
embrace the elementary abelian $2$-group of order $8$ and the group of
quaternions of order $8$. This is shown in \cite{PeVo}, section $8$. The
biggest group found so far, for which $\delta(G(\oa))\ne\delta_0(G(\oa))$ is
the cyclic group of order $21$ (see \cite{PeVo}, p.$36$).

Our goal is to prove that $\delta(G(\oa))=6p-18$ for each prime $p$ greater than $7$
(note that $\delta(G(\oa))\le 6p-18$ holds for each $p>7$). In order to achieve this
we need the following propositions:

\begin{lemma}\label{goodrow}
Suppose that $G(\oa)$, $G(\ob)$ are two groups of order $n$, and $a\oa b\ne a\ob b$
for some $a$, $b\in G$. Then $d(a)+d(b)+d(a\oa b)\ge n$.
\end{lemma}
\begin{proof}
\cite{PeVo} lemma $2.10$, or, more generally, \cite{DrM} lemma $2.4.$
\end{proof}

\begin{proposition}\label{drapal1}
Let $G(\oa)$, $G(\ob)$ be two groups.  Put $K=\{a\in G; d(a)<n/3\}$,
and assume that $|K|>3n/4$. Define a~mapping $f:G\to G$ by
$f(g)=a\ob b$ for any $g\in G$, $a$, $b\in K$, $g=a\oa b$. Then $f$ is
an isomorphism of $G(\oa)$ onto $G(\ob)$, and $f(a)=a$ for each $a\in K$.
Moreover, $f(g)\ne g$ for any $g\in G$ with $d(g)>2n/3$.
\end{proposition}
\begin{proof}
\cite{DrM} proposition $3.1.$
\end{proof}

\begin{proposition}\label{drapal2}
Let $G(\oa)$ be a~finite group of order $n\ge 5$. Then there exists a~transposition
$f$ of $G(\oa)$ with $m_f=\delta_0(G(\oa))$. Furthermore, $m_f\ge\delta_0(G(\oa))$
for any transposition $f$ of $G$. Finally, if $n\ge 12$, and $f$ is such
a~permutation of $G$ that $n>|\{g\in G; f(g)=g\}|>2n/3$, then
$m_f\ge\delta_0(G(\oa))$, and $f$ is a~transposition whenever
$m_f=\delta_0(G(\oa))$.
\end{proposition}
\begin{proof}
\cite{DrM} proposition $7.1.$
\end{proof}

\begin{lemma}\label{same1}
Assume that $G(\oa)$, $G(\ob)$ are two isomorphic groups of order $n>7$ satisfying
$dist(G(\oa),G(\ob))\le 6n-18$. Then we have $1_{G(\oa)}=1_{G(\ob)}$.
\end{lemma}
\begin{proof}
Let $e=1_{G(\oa)}$, $f=1_{G(\ob)}$. Assume that $e\ne f$. We would like
to prove that $d=dist(G(\oa),G(\ob))>6n-18$.

Put $E=\{(a,b)\in G\times G; \text{ }\{e,f\}\cap\{a,b\}\ne\emptyset\}$. We
show that $a\oa b\ne a\ob b$ for any $(a,b)\in E$. When $a=e$, we have $a\oa b=b$,
and $a\ob b\ne b$, since $a\ne f$. All remaining cases follow from symmetry.

For any $a\in G$ denote by $a^{-1}$, $a^{\ob}$ the inverse element of $a$ in
$G(\oa)$, $G(\ob)$, respectively. Define $I=\{a\in G; \text{ }a^{-1}=a^{\ob}\}$.

We prove that $d(a)\ge 4$ for any $a\in I$, $a\not\in\{e,f\}$.
Let $M=\langle e$, $f$, $a^{-1}$, $a^{-1}\oa f\rangle$ be an ordered set.
Note that all elements of $M$ are distinct. Hence also
$a\oa M=\langle a$, $a\oa f$, $e$, $f\rangle$ and
$a\ob M=\langle a\ob e$, $a$, $f$, $a\ob(a^{-1}\oa f)\rangle$ are
four-element sets. Moreover, each two respective elements of $a\oa M$ and $a\ob M$
are different.

If $a\not\in I$ and $b\in G$ are such that $a\oa b=a\ob b=c$, we have
$a^{\ob}\oa c\ne a^{\ob}\ob c$. Otherwise $b=a^{\ob}\ob a\ob b=a^{\ob}\ob c=
a^{\ob}\oa c\ne a^{-1}\oa c=b$, a contradiction. This means that
$d(a)+d(a^{\ob})\ge n$ for any $a\not\in I$.

Let $i=|I|$. We need to consider three possible cases.

$(i)$ Let $e\not\in I$, $f\not\in I$. If $i\ge n-4$, we have
$d\ge 4(n-4)+2n=6n-16>6n-18$. On the other hand, if $i\le n-5$,
then $d\ge (n-i)n/2+4i={n^2}/2+i(4-n/2)$. Since $n>7$, we can conclude
that $d\ge {n^2}/2+(n-5)(4-n/2)=13n/2-20>6n-18$.

$(ii)$ Let $|\{e,f\}\cap I|=1$. If $i\ge n-3$, then again (however,
the reason is different) $d\ge 4(n-4)+2n$. For $i\le n-4$, one
can see that $d\ge (n-i)n/2+4(i-1)+n={n^2}/2+i(4-n/2)-4+n
\ge{n^2}/2+(n-4)(4-n/2)-4+n=7n-20>6n-18$.

$(iii)$ Finally, let $\{e,f\}\subseteq I$. If $i\ge n-2$, we have
$d\ge 4(n-4)+2n$. If $i\le n-3$, then
$d\ge (n-i)n/2+4(i-2)+2n=n^2/2+i(4-n/2)-8+2n
\ge {n^2}/2+(n-3)(4-n/2)-8+2n=15n/2-20>6n-18$.

This proof can be found in \cite{PeVo}.
\end{proof}

Unfortunately, also some use of computers is needed in two special cases.

\section{Basic estimates}

From now on suppose that $G(\oa)$, $G(\ob)$ are two distinct groups of
prime order $p>7$. Let us denote by $H$ the set of all rows in
multiplication table of $G(\oa)$ at which
operations $\oa$ and $\ob$ completely agree, i.e. $H=\{g\in G; d(g)=0\}$.
Assume that $H$ is not empty, and $a$, $b$ belong to $H$. Then
$
(a\ob b)\oa g=(a\oa b)\oa g=a\oa (b\oa g)=a\oa (b\ob g)=a\ob (b\ob g)
=(a\ob b)\ob g=(a\oa b)\ob g$, which shows that $H$ is a~common subgroup
of $G(\oa)$ and $G(\ob)$.

According to lemma \ref{same1}, $H$ is never empty, when
$dist(G(\oa),G(\ob))<6p-18$.
Because there are no non-trivial subgroups in $\mathbb Z_p$, $H$ must be
the one element subgroup $1=1_{G(\oa)}=1_{G(\ob)}$, since
$G(\oa)$, $G(\ob)$ are distinct.

Put $m=min\{d(g); g\ne 1\}$. We know that $m>0$. The case $m=1$ is impossible,
hence $m>1$. In fact, as the following lemma shows, $m>2$.

\begin{lemma}\label{never2}
Let $G(\oa)$, $G(\ob)$ be two groups of odd order $n$. Then $d(g)\ne 2$
for any $g\in G$.
\end{lemma}
\begin{proof}
Let $\pi:G\to G$ be a~left translation by $g$ in $G(\oa)$, and
$\sigma:G\to G$ a~left translation by $g$ in $G(\ob)$. Then $g\oa a\ne g\ob a$
if and only if $\pi(a)\ne\sigma(a)$, i.e. ${\pi}^{-1}\circ\sigma(a)\ne a$.

Suppose that $d(g)=2$. This means that ${\pi}^{-1}\circ\sigma$ is a~transposition.
In particular, $sgn({\pi}^{-1}\circ\sigma)=-1$. But $sgn(\pi)={sgn(\pi)}^n=
sgn({\pi}^n)=sgn(id)=1$, and a~similar argument shows that also $sgn(\sigma)=1$,
a~contradiction.
\end{proof}

Suppose, for a~while, that $m\ge 6$. Then $dist(G(\oa),G(\ob))\ge 6(n-1)>6n-18$,
and we can see that this case is not interesting.

Some additional theory is needed for $m=3$, $4$, $5$.

We use symbol $\lceil x \rceil$ to denote the smallest integer $k$ such that
$x\le k$.

\begin{proposition}\label{estim1}
Let $G(\oa)$, $G(\ob)$ be two distinct groups of order $n\ge 5$. Then
either $dist(G(\oa),G(\ob))\ge \delta_0(G(\oa))$, or
\begin{displaymath}
    dist(G(\oa),G(\ob))\ge \lceil n/4\rceil\lceil n/3\rceil+
    (n-\lceil n/4\rceil-1)m.
\end{displaymath}
\end{proposition}

\begin{proof}
Put $K=\{a\in G; d(a)<n/3\}$.

$(i)$ Suppose that $|K|>3n/4$. By \ref{drapal1} there is an isomorphism
$f:G(\oa)\to G(\ob)$ such that $f(a)=a$ for each $a\in K$. If $n<12$, then we have
$|K|>3n/4>n-3$. Therefore $f$ must be a~transposition, and
$dist(G(\oa),G(\ob))=m_f\ge\delta_0(G(\oa))$ follows by \ref{drapal2}.
If $n\ge 12$, then
$dist(G(\oa),G(\ob))\ge\delta_0(G(\oa))$ follows at once from \ref{drapal2},
because $n>|K|>3n/4>2n/3$.

$(ii)$ Now, let $|K|\le 3n/4$. We show that there are at least
$\lceil n/4\rceil$ elements $g$ with $d(g)\ge\lceil n/3\rceil$.
Assume the contrary, i.e. assume that there are at least
$n-\lceil n/4\rceil+1$ elements $g$ with $d(g)<\lceil n/3\rceil$, so also with
$d(g)<n/3$. However, $n-\lceil n/4\rceil +1>3n/4$, a~contradiction
with $|K|\le 3n/4$.
\end{proof}

\begin{proposition}\label{estim2}
Let $G(\oa)$, $G(\ob)$ be as in previous proposition. Let's choose $h\in G$
such that $d(h)=m$, and $h_0,\dots,h_{m-1}$ are pairwise different elements
satisfying $h\oa h_i\ne h\ob h_i$ for $i=0,\dots,m-1$. Further suppose there is
an $l$-element subset $Y$ of $\{h_0,\dots,h_{m-1}\}$ such that
$Y\cap h\oa Y=\emptyset$. Then either $dist(G(\oa),G(\ob))\ge 6n-18$, or we get
\begin{align*}
dist(G(\oa),G(\ob))&\ge l(n-m)+(n-2l-1)m, \text{and} \tag 1\\
dist(G(\oa),G(\ob))&\ge l(n-m)+(\lceil n/4\rceil-2l)\lceil n/3 \rceil
+(n-\lceil n/4\rceil-1)m, \tag 2
\end{align*}
provided $\lceil n/4\rceil-2l\ge 0$.
\end{proposition}

\begin{proof}
Let us keep the notation of \ref{estim1}. If $|K|>3n/4$, then
$dist(G(\oa),G(\ob))\ge\delta_0(G(\oa))$ follows in the same way as in \ref{estim1}.
When $|K|\le 3n/4$, we have at least $\lceil n/4\rceil$ elements $g\in G$ for which
$d(g)\ge\lceil n/3\rceil$. Without loss of generality, put
$Y=\{h_0,\dots,h_{l-1}\}$. According to \ref{goodrow}, we get
\begin{align*}
d(h)+d(h_i)+d(h\oa h_i)\ge n, \text{or in other words}\\
d(h_i)+d(h\oa h_i)\ge n-m \text{ for each $i=0,\dots,l-1$}.
\end{align*}
This immediately proves \thetag{1}. In order to prove \thetag{2}, notice there
are at least $\lceil n/4\rceil -2l$ rows in $K$ not belonging to $Y\cup h\oa Y$.
\end{proof}

\begin{corollary}
When $G(\oa)$ is a~group of prime order $p>31$, then $\delta(G(\oa))=6p-18$.
\end{corollary}

\begin{proof}
Let $G(\ob)$ be the nearest group to $G(\oa)$. Since $m\ge 3$, it is easy
to see that we can always find a~set $Y$ (from \ref{estim2}) such that it has
at least two elements. Inequality \thetag{2} gives
\begin{displaymath}
    dist(G(\oa),G(\ob))\ge 2(p-m)+(\lceil p/4\rceil-4)\lceil p/3\rceil+
    (p-\lceil p/4\rceil-1)m.
\end{displaymath}
Observe that its right hand side is increasing in $m$. For $m=3$ we obtain
\begin{displaymath}
    dist(G(\oa),G(\ob))\ge 5p-9+(\lceil p/4\rceil-4)\lceil p/3\rceil-3\lceil p/4\rceil,
\end{displaymath}
and one can check that the expression on the r.h.s. is for $p>31$ always greater
than $6p-18$ (consider $p$ in form $12r+s$, say).
\end{proof}

\section{Case $m=5$}

Estimate \thetag{1} from \ref{estim2} turns out to be strong enough when $m=5$.
Let us denote,
for convenience, the powers of any $h$ in $G(\oa)$ by $h^r$. For example,
$h^2=h\oa h$.

\begin{lemma}
Let $G(\oa)$, $G(\ob)$ be two distinct groups of prime order $p>7$,
and suppose that $m=5$.
Then $dist(G(\oa),G(\ob))\ge 6p-18$.
\end{lemma}

\begin{proof}
Denote by $h$ one of the rows for which $d(h)=5$. Suppose that
$h^{i_0}$, $h^{i_1}$, $h^{i_2}$, $h^{i_3}$, $h^{i_4}$ are pairwise different
elements with $h\oa h^{i_j}\ne h\ob h^{i_j}$, $j=0,\dots,4$, where
$i_0<i_1<i_2<i_3<i_4<p$. We can suppose that $i_0>0$ (otherwise
$dist(G(\oa),G(\ob))\ge 6p-18$ follows from \ref{same1}).

We would like to find a~$3$-element subset $Y$ of
$\{h^{i_0}, h^{i_1}, h^{i_2}, h^{i_3}, h^{i_4}\}$ satisfying
$Y\cap h\oa Y=\emptyset$. Clearly, $h^{i_0+1}\ne h^{i_2}$, $h^{i_4}$. As
$i_0>0$, we have also $h^{i_2+1}$, $h^{i_4+1}\ne h^{i_0}$. Finally,
$h^{i_2+1}\ne h^{i_4}$, and $Y=\{h^{i_0}, h^{i_2}, h^{i_4}\}$ is such a~subset.
By \thetag{1} we know that
\begin{displaymath}
    dist(G(\oa),G(\ob))\ge 3(p-5)+(p-7)5=8p-50,
\end{displaymath}
and $8p-50$ is less than $6p-18$ only when $p<16$, i.e. $p\le 13$.

But when $p\le 13$ we have $dist(G(\oa),G(\ob))\ge 5p-5\ge 6p-18$.
\end{proof}

\section{Cases $m=4$, $m=3$}

\begin{proposition}\label{newprop}
For any two distinct groups $G(\oa)$, $G(\ob)$ of prime order $p>19$ with $m=4$
we have $dist(G(\oa),G(\ob))\ge 6p-18$.
\end{proposition}
\begin{proof}
Assume there is a~$3$-element subset $Y$ from \ref{estim2}. Then \thetag{1} yields
\begin{displaymath}
    dist(G(\oa),G(\ob))\ge3(p-4)+(p-7)4=7p-40,
\end{displaymath}
and $7p-40$ is less than $6p-18$ only when $p<22$, i.e. $p\le 19$. We cannot
improve this result by using estimate \thetag{2}, since $\lceil p/4\rceil\ge
2l=6$ if and only if $p\ge 21$.

It is not always feasible to find a~$3$-element subset $Y$ of
$\{h^{i_0},h^{i_1},h^{i_2},h^{i_3}\}$
with $Y\cap h\oa Y=\emptyset$.
One can show by tedious elementary methods
that this is not feasible if and only if $i_1=i_0+1$ and $i_3=i_2+1$. However, in
such a~case we can show that the transposition $f=(h^{i_1},h^{i_3})$ is
an isomorphism of $G(\oa)$ onto $G(\ob)$ (detailed proofs are given in \cite{PeVo}
$4.18$, $4.19$). Our wanted estimate then follows from \ref{drapal2}.
\end{proof}

There is no such estimate for $m=3$. We need more information about the group
operation $\ob$.

\begin{lemma}\label{struct}
Let $G(\oa)$, $G(\ob)$ be two groups of odd order $n$, and let $h$ be
a~common generator of $G(\oa)$, $G(\ob)$ with $d(h)=4$. Denote by
$h^{i_0}$, $h^{i_1}$, $h^{i_2}$, $h^{i_3}$ the pairwise different elements for which
$h\oa h^{i_j}\ne h\ob h^{i_j}$, $j=0,\dots,3$, where $i_0<i_1<i_2<i_3$. Then
$h\ob h^{i_0}=h\oa h^{i_2}$, $h\ob h^{i_2}=h\oa h^{i_0}$,
$h\ob h^{i_1}=h\oa h^{i_3}$, and $h\ob h^{i_3}=h\oa h^{i_1}$
\end{lemma}

\begin{proof}
Let $\pi$, $\sigma$ be as in the proof of \ref{never2}. Then ${\pi}^{-1}\circ\sigma$
is either a~$4$-cycle, or a~composition of two independent transpositions.
In fact, ${\pi}^{-1}\circ\sigma$ cannot be a~$4$-cycle, because
$sgn({\pi}^{-1}\circ\sigma)=1$. It is not difficult to observe that
${\pi}^{-1}\circ\sigma$ must be a~permutation $(i_0,i_2)(i_1,i_3)$.
\end{proof}

We can depict the situation as follows:
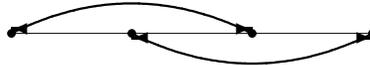
\begin{displaymath}
    \setlength{\unitlength}{1mm}\input{4arrows.lp}
\end{displaymath}
For $m=3$, the appropriate picture is (without proof):
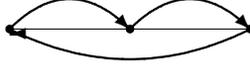
\begin{displaymath}
    \setlength{\unitlength}{1mm}\input{3arrows.lp}
\end{displaymath}

Now we have enough information to write efficient computer programs in order
to solve all remaining cases --- we only need to consider situations when
$m=4$ and $7<p<19$, or $m=3$ and $7<p<31$.

We will not give a~concrete implementation of requested algorithms (which
can be found in \cite{PeVo}), but we describe these algorithms in words instead.

Suppose that $p$ is a~prime between $7$ and $19$. We would like to modify the
canonical multiplication table of $\mathbb Z_p=G(\oa)$ in all possible ways,
such that the resulting table will be a~multiplication table of some group $G(\ob)$
satisfying $m=4$ (the other case $m=3$ is similar), and then check that
$dist(G(\oa),G(\ob))\ge 6p-18$.

By lemma \ref{same1}, the first row and the first column of $G(\oa)$ remain
unchanged. We choose some row $h\ne 0$ in $G$ and modify it at four places
$0<i_0<i_1<i_2<i_3<p$. According to \ref{struct}, this modification is given
by permutation $(i_0,i_2)(i_1,i_3)$, otherwise we never get a~group
multiplication table.

It is worth to point out that we do not need to go through all choices of $h\in G$.
In fact, we can fix only one row (a detailed explanation of this fact can be
found in \cite{PeVo}, $4.1$).
This trick speeds up the algorithm
$p-1$ times, and hence it is not essential.

Once we know one row of multiplication table of $G(\ob)$, we can build up $G(\ob)$
fully, because each non-zero element of $\mathbb Z_p$ is a~generator.

\section{Main result}

The algorithm described in section $5$ does not find any pair of groups $G(\oa)$,
$G(\ob)$ with $dist(G(\oa),G(\ob))<6p-18$, which, together with all previous results,
means that:

\begin{theorem}\label{main}
Each group of prime order $p>7$ has Cayley stability equal to $6p-18$.
\end{theorem}

Note that there are two groups $G(\oa)$, $G(\ob)$ of order $7$ with
$d(G(\oa),G(\ob))=18<24$ --- consider isomorphism
$f:G(\oa)\longrightarrow G(\ob)$ given by
\begin{displaymath}
    \begin{pmatrix}
        0&1&2&3&4&5&6\\
        0&1&4&5&2&3&6\\
    \end{pmatrix},
\end{displaymath}
so the estimate $p>7$ in \ref{main} cannot be improved.
These two groups are the nearest possible groups of order $7$ --- in other words,
$\delta(\mathbb Z_7)=18$.

It is easy to check that $\delta(\mathbb Z_2)=4$ and $\delta(\mathbb Z_3)=9$.
Computation reveals that $\delta(\mathbb Z_5)=12$. Here, the
group nearest to $\mathbb Z_5$ is obtained via
transposition $(2$, $3)$, for example.

\bibliographystyle{amsalpha}

\end{document}